\documentstyle[10pt]{amsart}
\newcommand{\Lam}{\Lambda}

\newcommand{\til}[1]{\tilde{#1}}
\newcommand{\innp}[2]{\left\langle #1,#2\right\rangle}

\newcommand{\im}[1]{\mbox{Im}\,{#1}}

\newcommand{\fr}[2]{\frac{#1}{#2}}
\newcommand{\alp}{\alpha}

\theoremstyle{definition}
\newtheorem{lemma}{Lemma}[section]

\newtheorem{theorem}[lemma]{Theorem}
\newtheorem{notation}[lemma]{Notation}
\newtheorem{remark}[lemma]{Remark}
\newtheorem{corollary}[lemma]{Corollary}

\newtheorem{definition}[lemma]{Definition}
\newcommand{\rarr}{\rightarrow}

\def\C{{\rm \kern.24em \vrule width .02em
   height1.4ex depth -.05ex \kern -.26em C}}

\newcommand{\ov}{\overline}

\newcommand{\noin}{\noindent}

\newcommand{\ext}{\bigwedge}

\newcommand{\itm}[1]{\item[{\rm (#1)}]}
\newcommand{\eucl}[2]{{\Bbb {#1}}^{#2}}
\newcommand{\project}[2]{{\Bbb P}_{\Bbb {#1}}^{#2}}

\newcommand{\spn}{{\Bbb C}\mbox{-span}}
\newcommand{\bun}[1]{{\cal O}_{#1}}
\begin{document}

\title{Unitary Tridiagonalisation in 
$M(4,{\Bbb C})$} 
\author{Vishwambhar Pati}
\address{Stat-Math Unit\\
Indian Statistical Institute\\
RVCE P.O., Bangalore 560 059, India}
\date{}
\begin{abstract} A question of interest in linear algebra is whether 
all $n\times n$ complex matrices can be unitarily 
tridiagonalised. The answer for all $n\neq 4$ (affirmative or negative) 
has been known for a while, whereas the case $n=4$ seems to have remained 
open. In this paper we settle the $n=4$ case in the affirmative. 
Some machinery from complex algebraic geometry needs to be used. 
\end{abstract}
\maketitle
\section{Main Theorem} Let $V=\eucl{C}{n}$, and $\innp{\;}{\;}$ be the usual 
euclidean hermitian inner product on $V$. $U(V)=U(n)$ denotes the group of 
unitary automorphisms of $V$ with respect to $\innp{\;}{\;}$.
$\{e_{i}\}_{i=1}^{n}$ will denote the standard orthonormal basis of $V$. 
$A\in \mbox{End}_{\Bbb C}(V)=
M(n,{\Bbb C})$ will always denote a ${\Bbb C}$-linear transformation of $V$.

A matrix $A=[a_{ij}]$ is said to be {\em tridiagonal} if $a_{ij}=0$ for all
$1\leq i,j\leq n$ such that $|i-j|\geq 2$. 

Then we have:

\begin{theorem} For $n\leq 4$, and $A\in M(n,{\Bbb C})$, there exists a
unitary $U\in U(n)$ such that $UAU^{\ast}$ is tridiagonal. 

\label{mainthm}
\end{theorem}

\bigskip

\begin{remark} The case $n=3$, and counterexamples for $n\geq 6$, are
 due to Longstaff, \cite{Long}. In the paper \cite{F-W}, Fong and Wu   
construct counterexamples for $n=5$, and provide a proof in certain special 
cases  for $n=4$. The article \S 4 of \cite{F-W} poses the $n=4$ case in general as an 
open question. Our main theorem above answers this question in the 
affirmative. In passing, we also provide another elementary proof for the 
$n=3$ case. 
\label{ngeq5rmk}
\end{remark}

\bigskip

\section{Some Lemmas} 
We need some preliminary lemmas, which we collect in this section. In
the sequel, we will also use the letter $A$ to denote the unique linear
transformation determined by $A$. 

\bigskip
\begin{lemma} Let $A\in M(n,{\Bbb C})$. For all $n$, the following are 
equivalent:

\bigskip

\begin{description}
\itm{i} There exists a unitary $U\in U(n)$ such that $UAU^{\ast}$ is
tridiagonal. 

\bigskip

\itm{ii} There exists a flag (=ascending sequence of ${\Bbb C}$-
subspaces) of $V=\eucl{C}{n}$:
$$
{0}=W_{0}\subset W_{1}\subset W_{2}\subset .......\subset W_{n}=V
$$
such that $\dim\,W_{i}=i,\; AW_{i}\subset W_{i+1}$ and
$A^{\ast}W_{i}\subset W_{i+1}$ for all $0\leq i\leq n-1$. 

\bigskip 

\itm{iii}There exists a flag in $V$:
$$
{0}=W_{0}\subset W_{1}\subset W_{2}\subset....\subset W_{n}=V
$$
such that $\dim\,W_{i}=i,\; AW_{i}\subset W_{i+1}$ and
$A(W_{i+1}^{\perp})\subset W_{i}^{\perp}$ for all $0\leq i\leq n-1$. 
\end{description}
\label{flaglemma}
\end{lemma}

\bigskip
\noin
{\bf Proof:} 

$(i)\Rightarrow (ii)$ 

Set $W_{i}=\spn(f_{1}, f_{2},...f_{i})$, where $f_{i}= U^{\ast}e_{i}$
and $e_{i}$ is the standard basis of $V=\eucl{C}{n}$. Since the 
matrix $[b_{ij}]:=UAU^{\ast}$
is tridiagonal, we have:
\begin{eqnarray*}
Af_{i}=b_{i-1,i}f_{i-1} + b_{ii}f_{i} + b_{i+1,i}f_{i+1}\;\;\;\;1\leq
i\leq n
\label{eqn1}
\end{eqnarray*}
(where $b_{ij}$ is understood to be $=0$ for $i,j\leq 0$ or $\geq n+1$).
Thus $AW_{i}\subset W_{i+1}$. Since $\{f_{i}\}_{i=1}^{n}$ is an
orthonormal basis for $V=\eucl{C}{n}$, we also have: 
\begin{eqnarray*}
A^{\ast}f_{i}=\ov{b}_{i,i-1}f_{i-1} + \ov{b}_{ii}f_{i} + \ov{b}_{i,i+1}f_{i+1}\;\;\;\;1\leq
i\leq n
\label{eqn2}
\end{eqnarray*}
which shows $A^{\ast}(W_{i})\subset W_{i+1}$ for all $i$ as well, and (ii) 
follows. 

\bigskip
$(ii)\Rightarrow (iii)$ 

$A^{\ast}W_{i}\subset
W_{i+1}$ implies $(A^{\ast}W_{i})^{\perp}\supset W_{i+1}^{\perp}$ for
$1\leq i\leq n-1$. But since $(A^{\ast}W_{i})^{\perp}=
A^{-1}(W_{i}^{\perp})$, we have $A(W_{i+1}^{\perp})\subset W_{i}^{\perp}$
for $1\leq i\leq n-1$ and (iii) follows. 

\bigskip
$(iii)\Rightarrow (i)$

Inductively choose an orthonormal basis $f_{i}$ of $V=\eucl{C}{n}$ so that 
$W_{i}$ is the span of $\{f_{1},..,f_{i}\}$. Since $A(W_{i})\subset W_{i+1}$, we
have:
\begin{eqnarray}
Af_{i}=a_{1i}f_{1}+a_{2i}f_{2}+....+a_{i+1,i}f_{i+1}
\label{Afi1}
\end{eqnarray}
Since $f_{i}\in (W_{i-1})^{\perp}$, and by hypothesis 
$A(W_{i-1}^{\perp})\subset W_{i-2}^{\perp}$, and 
$W_{i-2}^{\perp}=\spn(f_{i-1},f_{i},..,f_{n})$, we also have 
\begin{eqnarray}
Af_{i}=a_{i-1,i}f_{i-1}+a_{ii}f_{i}+...+a_{ni}f_{n}
\label{Afi2}
\end{eqnarray}
and by comparing the two equations (\ref{Afi1}), (\ref{Afi2}) above, it follows that 
$$
Af_{i}=a_{i-1,i}f_{i-1}+a_{ii}f_{i}+a_{i+1,i}f_{i+1}
$$
for all $i$, and defining the unitary $U$ by $U^{\ast}e_{i}=f_{i}$ makes
$UAU^{\ast}$ tridiagonal, so that (i) follows. \hfill $\Box$

\bigskip

\begin{lemma} Let $n\leq 4$. If there exists a $2$-dimensional
${\Bbb C}$-subspace $W$ of $V=\eucl{C}{n}$ such that $AW\subset W$ and
$A^{\ast}W\subset W$, then $A$ is unitarily tridiagonalisable. 
\label{invlemma}
\end{lemma}

\bigskip
\noin
{\bf Proof:} If $n\leq 2$, there is nothing to prove. For $n=3$ or $4$, 
the hypothesis implies that $A$ maps $W^{\perp}$ into
itself. Then, in an orthonormal basis $\{f_{i}\}_{i=1}^{n}$ of $V$ which
satisfies $W=\spn(f_{1},f_{2})$ and $W^{\perp}=\spn(f_{3},..,f_{n})$ the
matrix of $A$ is in $(1,2)$ (resp. $(2,2)$) block-diagonal form for
$n=3$ (resp. $n=4$), which is clearly tridiagonal. \hfill $\Box$ 

\bigskip

\begin{lemma}Every matrix $A\in M(3,{\Bbb C})$ is unitarily
tridiagonalisable.
\label{tridiag3}
\end{lemma}

\bigskip
\noin
{\bf Proof:}

For $A\in M(3,{\Bbb C})$, consider the homogeneous cubic
polynomial in $v=(v_{1},v_{2},v_{3})$ given by: 
$$
F(v_{1},v_{2},v_{3}):=\det(v,Av,A^{\ast}v)
$$
Note $v\wedge Av\wedge A^{\ast}v=
F(v_{1},v_{2},v_{3})e_{1}\wedge e_{2}\wedge e_{3}$. By a standard result 
in dimension theory (see \cite{Sha}, p. 74, Theorem 5) 
each irreducible component of 
$V(F)\subset \project{C}{2}$ is of dimension $\geq 1$, and $V(F)$ is 
non-empty. Choose some
$[v_{1}:v_{2}:v_{3}]\in V(F)$, and let $v=(v_{1},v_{2},v_{3})$ which is non-zero. 
Then we have the two cases:

\bigskip
\noin
{\em Case 1:}\;
$v$ is a common eigenvector for $A$ and $A^{\ast}$. Then the 2-dimensional 
subspace $W=({\Bbb C}v)^{\perp}$ is an invariant 
subspace for both $A$ and $A^{\ast}$, and applying the lemma \ref{invlemma} 
to $W$ yields the result. 

\bigskip
\noin
{\em Case 2:}\;
$v$ is not a common eigenvector for $A$ and $A^{\ast}$. Say it is not 
an eigenvector for $A$ (otherwise interchange the roles of $A$ and $A^{\ast}$). 
Set $W_{1}={\Bbb C}v,\;W_{2}=\spn(v, Av),W_{3}=V=\eucl{C}{3}$. Then $\dim\,W_{i}=i$,
for $i=1,2,3$, and the fact that $v\wedge Av\wedge A^{\ast}v=0$ shows
that $A^{\ast}W_{1}\subset W_{2}$. Thus, by (ii) of lemma \ref{flaglemma}, we
are done. 
\hfill $\Box$

Let $V=\eucl{C}{4}$ from now on, and $A\in
M(4,{\Bbb C})$. 

\begin{lemma} If $A$ and $A^{\ast}$ have a common eigenvector, then $A$
is unitarily tridiagonalisable.
\label{eigenlemma}
\end{lemma}

\bigskip
{\bf Proof:} 
If $v\neq 0$ is a common eigenvector for $A$ and $A^{\ast}$, the $3$-
dimensional subspace $W=({\Bbb C}v)^{\perp}$ is invariant under both $A$ and
$A^{\ast}$, and unitary tridiagonalisation of $A_{|W}$ exists from the
$n=3$ case of the lemma \ref{tridiag3} by a $U_{1}\in U(W)=U(3)$. The unitary 
$U=1\oplus U_{1}$ is the desired unitary in $U(4)$ tridiagonalising $A$. 
\hfill  $\Box$

\bigskip

\newpage

\begin{lemma} If the main theorem holds for all $A\in S$, where $S$ is 
{\em any} dense (in the classical topology) subset of 
$M(4, {\Bbb C})$, then it holds for all $A\in M(4,{\Bbb C})$.
\label{semisimp}
\end{lemma}

\bigskip
\noin
{\bf Proof:} 

This is a consequence of the compactness of the unitary group $U(4)$.
Indeed, let $T$ denote the closed subset of tridiagonal (with respect to the 
standard basis) matrices. 

\bigskip

Let $A\in M(4,{\Bbb C})$ be any general element. By the density of $S$,
there exist $A_{n}\in S$ such that
$A_{n}\rightarrow A$. By hypothesis, there are unitaries $U_{n}\in U(4)$ 
such that $U_{n}A_{n}U_{n}^{\ast}=T_{n}$, where $T_{n}$ are tridiagonal.
By the compactness of $U(4)$, and by passing to a subsequence if
necessary, we may assume that $U_{n}\rarr U\in U(4)$. Then
$U_{n}A_{n}U_{n}^{\ast}\rarr UAU^{\ast}$. That is $T_{n}\rarr
UAU^{\ast}$. Since $T$ is closed, and $T_{n}\in T$, we have $UAU^{\ast}$
is in $T$, viz., is tridiagonal. \hfill $\Box$

\bigskip

We shall now construct a suitable dense open subset $S\subset M(4, {\Bbb
C})$, and prove tridiagonalisability for a general $A\in S$ in the 
remainder of this paper. More precisely:

\bigskip

\begin{lemma} There is a dense open subset $S\subset M(4, {\Bbb C})$
such that:

\bigskip
\noin
\begin{description}
\itm{i} $A$ is nonsingular for all $A\in S$.

\bigskip

\itm{ii} $A$ has distinct eigenvalues for all $A\in S$.

\bigskip

\itm{iii} For each $A\in S$, the element $(t_{0}I + t_{1}A +
t_{2}A^{\ast})\in M(4,{\Bbb C})$ has rank $\geq 3$ for all 
$(t_{0},t_{1},t_{2})\neq (0,0,0)$ in $\eucl{C}{3}$. 

\end{description}
\label{genericset}
\end{lemma}

\bigskip
\noin
{\bf Proof:} The subset of singular matrices in $M(4,{\Bbb C})$ is
the complex algebraic subvariety of complex codimension one defined by 
$Z_{1}=\{A:\det\,A=0\}$. 
Let $S_{1}$, (which is just $GL(4,{\Bbb C})$) be its complement. Clearly $S_{1}$ is
open and dense in the classical topology (in fact, also in the Zariski 
topology). 

\bigskip

A matrix $A$ has distinct eigenvalues iff its characteristic polynomial 
$\phi_{A}$ has distinct roots. This happens iff the discriminant
polynomial of $\phi_{A}$, which is a 4-th degree homogeneous polynomial
$\Delta(A)$ in the entries of $A$, is not zero. The zero set $Z_{2}=V(\Delta)$ is 
again a codimension-1 subvariety in $M(4,{\Bbb C})$, so its complement 
$S_{2}=\left(V(\Delta)\right)^{c}$ is open and dense in both the classical and
Zariski topologies. 

\bigskip 

To enforce (iii), we claim that the set defined by 
$$
Z_{3}:=\{A\in M(4,{\Bbb C}):\,\mbox{rank}\,(t_{0}I + t_{1}A + t_{2}A^{\ast})
\leq 2\;\mbox{for some}\;(t_{0},t_{1},t_{2})\neq (0,0,0)\;\;\mbox{in}\;\;
\eucl{C}{3}\} 
$$
is a proper {\em real} algebraic subset of $M(4,{\Bbb C})$. The proof hinges 
on the fact that three general cubic curves in $\project{C}{2}$ having a 
point in common imposes an algebraic condition on their coefficients. 

\bigskip

Indeed, saying that $\mbox{rank}\,(t_{0}I + t_{1}A + t_{2}A^{\ast})
\leq 2\;\mbox{for some}\;(t_{0},t_{1},t_{2})\neq (0,0,0)$ is equivalent
to saying that the third exterior power 
$\ext^{3}(t_{0}I + t_{1}A + t_{2}A^{\ast})$ is the zero map, for some 
$(t_{0},t_{1},t_{2})\neq 0$. 
This is equivalent to demanding that there exist a 
$(t_{0},t_{1},t_{2})\neq 0$ such that the determinants of 
all the $3\times 3$-minors of $(t_{0}I + t_{1}A + t_{2}A^{\ast})$ are zero. 

\bigskip

Note that the (determinants of) the $(3\times 3)$-minors of 
$(t_{0}I + t_{1}A + t_{2}A^{\ast})$,
denoted as $M_{ij}(A,t)$ (where the $i$-th row and $j$-th column are deleted) 
are complex valued, complex algebraic and ${\Bbb C}$-homogeneous of degree 3 in 
$t=(t_{0},t_{1},t_{2})$, with coefficients 
real algebraic of degree 3 in the variables $(A_{ij},\bar{A}_{ij})$ (or, equivalently, in
$\mbox{Re}\,A_{ij}, \mbox{Im}\,A_{ij}$), where $A=[A_{ij}]$.   

\bigskip

We know that the space of all homogeneous polynomials of degree 3 with complex 
coefficients in $(t_{0},t_{1},t_{2})$ (upto scaling) is parametrised by the projective
space $\project{C}{9}$ (the Veronese variety, see \cite{Sha}, p.52). 
We first consider the complex algebraic variety:

$$
X=\{(P,Q,R,[t])\in \project{C}{9}\times 
\project{C}{9}\times \project{C}{9}\times \project{C}{2}:
P(t)=Q(t)=R(t)=0\}
$$
where $[t]:=[t_{0}:t_{1}:t_{2}]$, and $(P, Q, R)$ denotes a triple of
homogeneous polynomials. This is just the subset of those 
$(P,Q,R,[t])$ in the
product $\project{C}{9}\times 
\project{C}{9}\times \project{C}{9}\times \project{C}{2}$ such that the 
point $[t]$ lies on all three of the plane cubic curves $V(P),V(Q),V(R)$. 
Since $X$ is defined by bihomogenous degree (1,1,1,3) equations, it is a complex
algebraic subvariety of the quadruple product. Its image under the first projection 
$Y:=\pi_{1}(X)\subset \project{C}{9}\times \project{C}{9}\times 
\project{C}{9}$ is 
therefore an algebraic subvariety inside this triple product (see 
\cite{Sha}, p. 58, Theorem 3). $Y$ is a
proper subvariety because, for example, the cubic polynomials
$P=t_{0}^{3}, Q=t_{1}^{3}, R=t_{2}^{3}$ have no common non-zero root. 

\bigskip

Denote pairs $(i,j)$ with $1\leq i,j\leq 4$ by capital letters like $I,J,K$ etc. 
From the minorial determinants $M_{I}(A,t)$, we can define various 
{\em real algebraic} maps:

\begin{eqnarray*}{}
\Theta_{IJK}:M(4,{\Bbb C})&\rarr & \project{C}{9}\times 
\project{C}{9}\times \project{C}{9}\\
A &\mapsto & (M_{I}(A,t),M_{J}(A,t),M_{K}(A,t))
\end{eqnarray*}
for $I,J,K$ distinct. 
Clearly, $\ext^{3}(t_{0}I+t_{1}A+t_{2}A^{\ast})=0$ for some 
$t=(t_{0},t_{1},t_{2})\neq (0,0,0)$ iff $\Theta_{IJK}(A)$ lies in   
in the complex algebraic subvariety $Y$ of $\project{C}{9}\times 
\project{C}{9}\times \project{C}{9}$, for all $I,J,K$ distinct. Hence   
the subset $Z_{3}\subset M(4,{\Bbb C})$ defined above is the intersection:

$$
Z_{3}=\bigcap_{I,J,K}\Theta_{IJK}^{-1}(Y)
$$
where $I,J,K$ runs over all distinct triples of pairs $(i,j),\; 1\leq
i,j\leq 4$. 

\bigskip

We claim that $Z_{3}$ is a proper real algebraic subset of 
$M(4,{\Bbb C})$. Clearly, since each $M_{I}(A,t)$ is real algebraic 
in the variables $\mbox{Re}\,A_{ij},\,\mbox{Im}\,A_{ij}$ the map $\Theta_{IJK}$ is real
algebraic. Since $Y$ is complex and hence real algebraic, its
inverse image $\Theta_{IJK}^{-1}(Y)$, defined by the real algebraic 
equations obtained upon substitution of the
components $M_{I}(A,t),M_{J}(A,t),M_{K}(A,t)$ in the equations that
define $Y$, is also real algebraic. Hence the set $Z_{3}$ is a 
real algebraic subset of $M(4,{\Bbb C})$.

\bigskip

To see that $Z_{3}$ is a {\em proper} subset of $M(4, {\Bbb C})$, we simply
consider the matrix (defined with respect to the standard orthonormal 
basis $\{e_{i}\}_{i=1}^{4}$ of $\eucl{C}{4}$):

$$
A=
\left[
\begin{array}{cccc}
0& 1& 0& 0\\
0& 0& 1& 0\\
0& 0& 0& 1\\
0& 0& 0& 0
\end{array}
\right]
$$ 
For $t=(t_{0},t_{1},t_{2})\neq 0$, we see that:
$$
t_{0}I + t_{1}A + t_{2}A^{\ast}= 
\left[
\begin{array}{cccc}
t_{0}& t_{1}& 0& 0\\
t_{2}& t_{0}& t_{1}& 0\\
0& t_{2}& t_{0}& t_{1}\\
0& 0& t_{2}& t_{0}
\end{array}
\right]
$$
For this matrix above, the minorial determinant $M_{41}(A,t)=
t_{1}^{3}$, whereas $M_{14}(A,t)=t_{2}^{3}$. The only common
zeros to these two minorial determinants are points $[t_{0}:0:0]$. Setting $t_{1}=t_{2}=0$ in the
matrix above gives $M_{ii}(A,t)=t_{0}^{3}$ for $1\leq i\leq 4$. 
Thus $t_{0}$ must also be $0$ for all the minorial determinants to vanish. 
Hence the matrix $A$ above lies outside the 
real algebraic set $Z_{3}$.

\bigskip

It is well known that a proper real algebraic subset in euclidean
space cannot have a non-empty interior. Thus the complement $Z_{3}^{c}$ 
is dense and open in the classical and real-Zariski topologies. 
Take $S_{3}=Z_{3}^{c}$. 

\bigskip

Finally, set 
$$
S:=S_{1}\cap S_{2}\cap S_{3}=\left(\bigcup_{i=1}^{3} Z_{i}\right)^{c}
$$
which is also open and
dense in the classical topology in $M(4, {\Bbb C})$. Hence the lemma.
$\Box$

\bigskip

\begin{remark} One should note here that for {\em each} matrix 
$A\in M(4,{\Bbb C})$, there will be at least a curve of points
$[t]=[t_{0}:t_{1}:t_{2}]\in \project{C}{2}$ (defined by the vanishing
of $\det\,(t_{0}I+t_{1}A+t_{2}A^{\ast})$), on which  
$(t_{0}I+t_{1}A+t_{2}A^{\ast})$ is singular. Similarly for each $A$ there is 
at least a curve of points on which the trace $\mbox{tr}\left(
\ext^{3}(t_{0}I+t_{1}A+t_{2}A^{\ast})\right)$
vanishes, and so a non-empty (and generally a finite) set on which 
{\em both} these polynomials vanish, by dimension theory (\cite{Sha}, Theorem 5,
p. 74). Thus for {\em each} $A\in M(4, {\Bbb C})$, there is a 
at least a non-empty finite set of points 
$[t]$ such that $(t_{0}I+t_{1}A+t_{2}A^{\ast})$ has $0$ as a repeated
eigenvalue. For example, for the matrix $A$ constructed at the end of the 
previous lemma, we see that the matrix 
$(t_{0}I +t_{1}A+t_{2}A^{\ast})$ is strictly upper-triangular and thus has $0$ as an 
eigenvalue of multiplicity $4$ for all $(0,t_{1},0)\neq 0$, but nevertheless has rank $3$ for all 
$(t_{0},t_{1},t_{2})\neq (0,0,0)$. 

\bigskip

Indeed, as (iii) of the lemma above shows, for $A$ in the open 
dense subset $S$, the kernel $\ker\,(t_{0}I+t_{1}A+t_{2}A^{\ast})$ is at most 
1-dimensional {\em for all} $[t]=[t_{0}:t_{1}:t_{2}]\in \project{C}{2}$. 
\label{genericityremk}
\end{remark}

\bigskip

\section{The varieties $C,\;\Gamma,\;$ and $D$}
\begin{notation} 

In the light of the lemmas \ref{semisimp} and 
\ref{genericset} above, we shall henceforth assume $A\in S$. 
As is easily verified, this implies $A^{\ast}\in S$ as 
well. We will also henceforth assume, in view of lemma ~\ref{eigenlemma}
above, that {\em $A$ and $A^{\ast}$ have no common eigenvectors}. (For example, 
this rules out $A$ being normal, in which case we know that 
the main result for $A$ is true by the spectral theorem). Also, in view
of lemma \ref{invlemma}, we shall assume that $A$ and $A^{\ast}$ do not
have a common 2-dimensional invariant subspace. 

\bigskip

In $\project{C}{3}$, the complex projective space of $V=\eucl{C}{4}$, 
we denote the equivalence class of $v\in V\setminus 0$ by $[v]$. For a 
$[v]\in \project{C}{3}$, we define $W([v])$ (or simply $W(v)$ when no confusion
is likely) by: 
$$
W([v]):=\spn (v, Av, A^{\ast}v)
$$
Since we are assuming that $A$ and $A^{\ast}$ have no common
eigenvectors, we have $\dim\,W([v])\geq 2$ for all $[v]\in
\project{C}{3}$.

\bigskip

Denote the 
four distinct points in $\project{C}{3}$ representing the
four distinct eigenvectors of $A$ (resp. $A^{\ast}$) by $E$ (resp.
$E^{\ast}$). By our assumption above, $E\cap E^{\ast}=\phi$. 
\label{assump}
\end{notation}

\bigskip

\begin{lemma} Let $A\in M(4,{\Bbb C})$ be as in \ref{assump} above. Then the closed
subset:
$$
C=\{[v]\in\project{C}{3}: 
v\wedge Av\wedge A^{\ast}v=0\}
$$
is a closed projective variety. This variety $C$ is precisely the subset of $[v]\in\project{C}{3}$ for 
which the dimension $\dim W([v])=\dim\left(\spn\,\{v,Av,A^{\ast}v\}\right)$ is 
exactly 2. 
\label{2lemma}
\end{lemma}

\bigskip
\noin
{\bf Proof:} 
That $C$ is a closed projective variety is clear from the fact that it is
defined as the set of common zeros of all the four $(3\times 3)$-minorial
determinants of the $(3\times 4)$-matrix 
$$
\Lambda:=
\left[
\begin{array}{c}
v\\
Av\\
A^{\ast}v
\end{array}
\right]
$$ 
(which are all degree-$3$ homogeneous polynomials in the components
of $v$ with respect to some basis). Also $C$ is nonempty since 
it contains $E\cup E^{\ast}$. 

\bigskip

Also, since $A$ and $A^{\ast}$ are nonsingular by the assumptions 
in \ref{assump}, the wedge 
product $v\wedge Av\wedge A^{\ast}v$ of the three non-zero vectors $v,
Av, A^{\ast}v$ vanishes precisely when the space $W([v])=\spn\{v, Av,
A^{\ast}v\}$ is of dimension $\leq 2$. Since by \ref{assump}, $A$,
$A^{\ast}$ have no common eigenvectors, the dimension $\dim\,W([v])\geq 2$
for all $[v]\in \project{C}{3}$, so $C$ is precisely the locus of
$[v]\in \project{C}{3}$ for which the space $W([v])$ is 2-dimensional.
\hfill $\Box$

\bigskip

Now we shall show that for $A$ as in \ref{assump}, the variety $C$ defined
above is of pure dimension one. For this, we need to define some more
associated algebraic varieties and regular maps. 

\bigskip

\begin{definition} Let us define the bilinear map:
\begin{eqnarray*}
B:\eucl{C}{4}\times \eucl{C}{3} &\rarr& \eucl{C}{4}\\
(v, t_{0},t_{1}, t_{2}) &\mapsto& B(v,t):=(t_{0}I+t_{1}A+t_{2}A^{\ast})v
\label{bilinear}
\end{eqnarray*}
We then have the linear maps $B(v,-):\eucl{C}{3}\rarr \eucl{C}{4}$ for 
$v\in \eucl{C}{4}$ and $B(-,t):\eucl{C}{4}\rarr \eucl{C}{3}$ for $t\in 
\eucl{C}{3}$. 

\bigskip

Note that the image $\im\,B(v,-)$ is the span of $\{v, Av,
A^{\ast}v\}$, which was defined to be $W(v)$. For a fixed $t$, denote the
kernel 
$$
K(t):=\ker(B(-,t):\eucl{C}{4}\rarr\eucl{C}{4})
$$  

\bigskip
\noin

Denoting $[t_{0}:t_{1}:t_{2}]$ by $[t]$ and 
$[v_{1}:v_{2}:v_{3}:v_{4}]$ by $[v]$ for brevity, we define:

$$
\Gamma :=\{([v],[t])\in \project{C}{3}\times
\project{C}{2}: B(v,t)=0\}
$$

\bigskip

Finally, define the variety $D$ by:

$$
D\subset \project{C}{2}:=\{[t]\in \project{C}{2}:\; \det\,B(-,t)=\det\,(t_{0}I+t_{1}A+t_{2}A^{\ast})=0\}
$$

\bigskip
\noin
Let 
$$
\pi_{1}:\project{C}{3}\times \project{C}{2}\rarr \project{C}{3},\;\;\;
\pi_{2}:\project{C}{3}\times \project{C}{2}\rarr \project{C}{2}
$$
denote the two projections. 
\label{defgamma}
\end{definition} 

\bigskip

\begin{lemma}

We have the following facts:

\bigskip
\noin
\begin{description}
\itm{i} $\pi_{1}(\Gamma)=C$, and $\pi_{2}(\Gamma)=D$. 
\bigskip
\itm{ii} $\pi_{1}:\Gamma\rarr C$ is 1-1, and the map $g$ defined by
$$
g:=\pi_{2}\circ\pi_{1}^{-1}:C\rarr D   
$$
is a regular map so that $\Gamma$ is the graph of $g$ and isomorphic as a variety 
to $C$.

\bigskip
\noin
\itm{iii} $D\subset \project{C}{2}$  is a plane curve, 
of pure dimension one. The map $\pi_{2}:\Gamma\rarr D$ is 1-1, and the
map $\pi_{1}\circ \pi_{2}^{-1}: D\rarr C$ is the regular inverse of the 
regular map $g$ defined above in (ii). Again $\Gamma$ is also the graph of 
this regular inverse $g^{-1}$, and $D$ and $\Gamma$ are isomorphic as varieties.
In particular, $C$ and $D$ are isomorphic as varieties, 
and thus $C$ is a curve in $\project{C}{3}$ 
of pure dimension one. 

\bigskip
\noin
\itm{iv} Inside 
$\project{C}{3}\times \project{C}{2}$, each irreducible component of the 
intersection of the four divisors $D_{i}:=(B_{i}(v,t)=0)$ for $i=1,2,3,4$  (where $B_{i}(v,t)$ is the 
$i$-th component of $B(v,t)$ with respect to a fixed basis of 
$\eucl{C}{4}$) occurs
with multiplicity 1. (Note that $\Gamma$ is set-theoretically the
intersection of these four divisors, by definition). 

\end{description}
\label{poslemma2}
\end{lemma}

\bigskip
\noin
{\bf Proof:} 

It is clear that $\pi_{1}(\Gamma)=C$, because 
$B(v,t)=t_{0}v+t_{1}Av+t_{2}A^{\ast}v=0$ for some $[t_{0}:t_{1}:t_{2}]\in 
\project{C}{2}$ iff $\dim\,W(v)\leq 2$, and since 
$A$ and $A^{\ast}$ have no common eigenvectors, this means $\dim\,W(v)=2$. 
That is, $[v]\in C$. 

\bigskip

Clearly $[t]\in \pi_{2}(\Gamma)$ iff there exists a 
$[v]\in \project{C}{3}$ such that $B(v,t)=0$. That is, iff 
$\dim\,\ker\,B(-,t)\,\geq 1$, that is, iff 
$$
G(t_{0},t_{1},t_{2}):=\det\,B(-,t)=0
$$
Thus $D=\pi_{2}(\Gamma)$ and is defined by a single degree 
4 homogeneous polynomial $G$ inside $\project{C}{2}$. It is a curve
of pure dimension 1 in $\project{C}{2}$ by standard dimension theory 
(see \cite{Sha}, p. 74, Theorem 5) because, for 
example $[1:0:0]\not\in D$ so $D\neq \project{C}{2}$. So 
$\pi_{2}(\Gamma)=D$, and this proves (i).

\bigskip

To see (ii), for a given $[v]\in C$, we claim there is exactly one 
$[t]$ such that $([v],[t])\in \Gamma$. Note that $([v],[t])\in \Gamma$ iff 
the linear map:
\begin{eqnarray*}
B(v,-):\eucl{C}{3} &\rarr & \eucl{C}{4}\\
t &\mapsto& (t_{0}I +t_{1}A +t_{2}A^{\ast})v
\end{eqnarray*}
has a non-trivial kernel containing the line ${\Bbb C}t$. 
That is, $\dim\,\im\,B(v,-)\leq 2$. But the image $\im\,B(v,-)=W(v)$, 
which is of dimension $2$ for all $v\in C$ by our assumptions. 
Thus its kernel must be exactly one dimensional, defined 
by $\ker\,B(v,-)={\Bbb C}t$. Thus $([v],[t])$ is the unique point in 
$\Gamma$ lying in $\pi_{1}^{-1}[v]$, viz. for each $[v]\in C$, the vertical line 
$[v]\times \project{C}{2}$ intersects $\Gamma$ in a single point, call
it $([v], g[v])$. So $\pi_{1}:\Gamma\rarr C$ is 1-1, and $\Gamma$ is the graph of a 
map $g:C\rarr D$. Since $g([v])=\pi_{2}\pi_{1}^{-1}([v])$ for 
$[v]\in C$, and $\Gamma$ is algebraic, $g$ is a regular map. This proves 
(ii).

\bigskip

To see (iii), note that for $[t]\in D$, by definition, 
the dimension $\dim\,\ker\,B(-,t)\geq 1$. By the fact that $A\in S$, and 
(iii) of the lemma \ref{genericset}, we know that $\dim\,\ker\,B(-
,t)\leq 1$ for all $[t]\in \project{C}{2}$. Thus, denoting 
$K(t):=\ker\,B(-,t)$ for $[t]\in D$, we have:
\begin{eqnarray}
\dim K(t)=1 \;\;\;\;\mbox{for all}\;\;\;\;t\in D
\label{kerdimension}
\end{eqnarray}

\bigskip

Hence we see that the unique 
projective line $[v]$ corresponding to ${\Bbb C}v = K(t)$ yields the 
unique  element of $C$, such that $([v],[t])\in \Gamma$. Thus
$\pi_{2}:\Gamma\rarr D$ is 1-1, and the regular map
$\pi_{1}\circ\pi_{2}^{-1}:D\rarr C$ is the regular inverse to the map $g$ of (ii)
above. $\Gamma$ is thus also the graph of $g^{-1}$ and, in particular, 
is isomorphic to $D$. Since $g$ is an isomorphism of curves, and $D$ is of 
pure dimension 1, it follows that $C$ is of pure dimension one. This proves
(iii). 

\bigskip

To see (iv), we need some more notation. 

Note that 
$D\subset \project{C}{2}\setminus \{[1;0;0]\}$, 
(because there exists no $[v]\in \project{C}{3}$ such that $I.v=0$!). 
Thus there is a regular map:
\begin{eqnarray}
\theta: D &\rarr& \project{C}{1}\\ \nonumber
[t_{0}:t_{1}:t_{2}] &\mapsto& [t_{1}:t_{2}]
\end{eqnarray}
Let $\Delta(t_{1},t_{2})$ be the discriminant polynomial of the
characteristic polynomial $\phi_{t_{1}A+t_{2}A^{\ast}}$ of 
$t_{1}A+t_{2}A^{\ast}$. Clearly
$\Delta(t_{1},t_{2})$ is a homogeneous polynomial of degree $4$ in
$(t_{1},t_{2})$, and it is not the zero polynomial because, for example, 
 $\Delta(1,0)\neq 0$, for $\Delta(1,0)$ is the discriminant of 
$\phi_{A}$, which has distinct roots (=the distinct eigenvalues of $A$) 
by the assumptions \ref{assump} on $A$. Let $\Sigma\subset \project{C}{1}$ be the
zero locus of $\Delta$, which is a finite set of points.  Note that the 
fibre $\theta^{-1}([1:\mu])$ consists of all $[t:1:\mu]\in D$ such that 
$-t$ is an eigenvalue of $A+\mu A^{\ast}$, which are at most four in number. Similarly 
the fibres $\theta^{-1}([\lambda:1])$ are also finite. Thus the subset
of $D$ defined by:

$$
F:=\theta^{-1}(\Sigma)
$$
is a finite subset of $D$. $F$ is precisely the set of points
$[t]=[t_{0}:t_{1}:t_{2}]$ such that 
$B(-,t)=(t_{0}I+t_{1}A+t_{2}A^{\ast})$ has $0$ as a repeated eigenvalue. 

\bigskip

Since $\pi_{2}:\Gamma\rarr D$ is 1-1, the inverse image:
$$
F_{1}=\pi_{2}^{-1}(F) \subset \Gamma
$$ 
is a finite subset of $\Gamma$. 

\bigskip

We will now prove that for each
irreducible component $\Gamma_{\alp}$ of $\Gamma$, and each 
point $x=([a],[b])$ in 
$\Gamma_{\alp}\setminus F_{1}$, the four equations 
$\{B_{i}(v,t)=0\}_{i=1}^{4}$ are the generators of the 
ideal of the variety $\Gamma_{\alp}$ in an affine neighbourhood of $x$, 
where $B_{i}(v,t)$ are the components of $B(v,t)$ with respect to a fixed 
basis of $\eucl{C}{4}$. Since $F_{1}$ is a finite set, this
will prove (iv), because the multiplicity of $\Gamma_{\alp}$ in the
intersection cycle of the four divisors $D_{i}=(B_{i}(v,t)=0)$ in
$\project{C}{3}\times \project{C}{2}$ is determined by generic points
on $\Gamma_{\alp}$, for example all points of $\Gamma_{\alp}\setminus F_{1}$. We will
prove this by showing that for 
$x=([a],[b])\in \Gamma_{\alp}\setminus F_{1}$, the four divisors
$(B_{i}(v,t)=0)$ intersect transversely at $x$. 

\bigskip

So let $\Gamma_{\alp}$ be some irreducible component of $\Gamma$, with 
$x=([a],[b])\in \Gamma_{\alp}\setminus F_{1}$. 

\bigskip

Fix an $a\in \eucl{C}{4}$ representing 
$[a]\in C_{\alp}:=\pi_{1}(\Gamma_{\alp})$, and also fix 
$b\in \eucl{C}{3}$ representing $[b]=g([a])\in g(C_{\alp})$. Also fix a 
3-dimensional linear complement
$V_{1}:=T_{[a]}(\project{C}{3})\subset \eucl{C}{4}$ to $a$ and similarly, 
fix a 2-dimensional linear complement $V_{2}=T_{[b]}(\project{C}{2})\subset \eucl{C}{3}$
to $b$. (The notation comes from the fact that $T_{[v]}(\project{C}{n})\simeq 
\eucl{C}{n+1}/{\Bbb C}v$, which we are identifying non-canonically with these respective 
complements $V_{i}$). These complements also provide local coordinates in the
respective projective spaces as follows. Set coordinate charts $\phi$   
around $[a]\in\project{C}{3}$ by $[v]=\phi(u):=[a+u]$, and $\psi$ around 
$[b]\in\project{C}{2}$ by $[t]=\psi(s):=[b+s]$, 
where $u\in V_{1}\simeq \eucl{C}{3}$, and $s\in V_{2}\simeq \eucl{C}{2}$. 
The images $\phi(V_{1})$ and $\psi(V_{2})$ are affine neighbourhoods of 
$[a]$ and $[b]$ respectively. These charts are 
like `stereographic projection' onto the tangent space and depend on the 
initial choice of $a$ (resp. $b$) representing $[a]$ (resp. $[b]$), 
and are {\em not} the standard coordinate systems on projective space, but 
more convenient for our purposes. 

\bigskip 

Then the local affine representation of $B(v,t)$ on the affine open 
$V_{1}\times V_{2}=\eucl{C}{3}\times\eucl{C}{2}$, which we denote by 
$\beta$, is given by:
$$
\beta(u,s):=B(a+u, b+s)
$$

\bigskip

Note that $\ker\,B(a,-)={\Bbb C}b$, where $[b]=g([a])$, so that $B(a,-)$
passes to the quotient as an isomorphism:

\begin{eqnarray}
B(a,-):V_{2}\;\widetilde{\longrightarrow}\; W([a])
\label{imageBa}
\end{eqnarray}
where $W([a])$ is 2-dimensional. 

\bigskip

Similarly, since $B(-,b)$ has one dimensional kernel ${\Bbb C}a=K(b)\subset
\eucl{C}{4}$, by (\ref{kerdimension}) above, 
we also have the other isomorphism:
\begin{eqnarray}
B(-,b):V_{1}\;\widetilde{\longrightarrow}\; \im\,B(-,b)
\label{imageBb}
\end{eqnarray}
where $\im B(-,b)$ is 3-dimensional, therefore. 

\bigskip

Now one can easily calculate the derivative $D\beta(0,0)$ of $\beta$ at 
$(u,s)=(0,0)$. Let $(X,Y)\in V_{1}\times V_{2}$. Then, by bilinearity of 
$B$, we have: 
\begin{eqnarray*}
\beta(X,Y)-\beta(0,0)&=&B(a+X,b+Y)-B(a,b)\\
&=& B(X,b) + B(a,Y) + B(X,Y)
\label{betaeqn}
\end{eqnarray*}

\bigskip

Now since $B(X,Y)$ is quadratic, it follows that:
\begin{eqnarray}
D\beta(0,0):V_{1}\times V_{2}&\rarr & \eucl{C}{4}\nonumber\\
(X,Y) &\mapsto& B(X,b)+B(a,Y)
\label{derivbeta}
\end{eqnarray}

\bigskip

By the equations (\ref{imageBa}) and (\ref{imageBb}) above, we see that
the image of $D\beta(0,0)$ is precisely $\im B(-,b)+W([a])$. 

\bigskip
\noin
{\bf Claim:} For $([a],[b])\in \Gamma_{\alp}\setminus F_{1}$, the space 
$\im B(-,b)+W([a])$ is all of $\eucl{C}{4}$. 

\bigskip
\noin
{\bf Proof of Claim:}
Denote $T:=B(-,b)$ for brevity.
Clearly $a\in W([a])$ by definition of $W([a])$. Also, $a\in \ker T=K(b)$. 
 Then it is enough to show that $0\neq a$ is not 
in the image of $T$. For, if 
 $a\in \im T$, we would have
$a=Tw$ for some $w\not\in K(b)=\ker T$ and $w\neq 0$. In fact $w$ is not a 
multiple of $a$ since $Tw=a\neq 0$ whereas $a\in \ker T$. Thus we would have $T^{2}w=0$,
and completing $e_{1}=a=Tw, e_{2}=w$ to a basis $\{e_{i}\}_{i=1}^{4}$ 
of $\eucl{C}{4}$, the matrix of $T$
with respect to such a basis would be of the form:
$$
\left[
\begin{array}{cccc}
0& 1& *& *\\
0& 0& *& *\\
0& 0& *&  *\\
0& 0& *&  *
\end{array}
\right]
$$

Thus $T=B(-,b)$ would have $0$ as a repeated eigenvalue. But we have
stipulated that $([a],[b])\not\in F_{1}$, so that $[b]\not\in F$, and
hence $B(-,b)$ does not have $0$ as a repeated eigenvalue. Hence 
the non-zero vector $a\in W([a])$ is not in $\im T$. Since $\im T$ is 3-dimensional, we have
$\eucl{C}{4}=\im T +W([a])$, and this proves the claim.
\hfill
$\Box$

\bigskip

In conclusion, all the points of $\Gamma_{\alp}\setminus F_{1}$ are in
fact smooth points of $\Gamma_{\alp}$, and the local equations for
$\Gamma_{\alp}$ in a small neighbourhood of such a point are precisely 
the four equations $\beta_{i}(u,s)=0$, $1\leq i \leq 4$. This proves (iv), and the lemma.
\hfill $\Box$

\section{Some algebraic bundles}

We construct an algebraic line bundle with 
a (regular) global section over $C$. By showing that 
this line bundle has positive degree, we will conclude that the section 
has zeroes in $C$. Any zero of this
section will yield a flag of the kind required by lemma \ref{flaglemma}. 
One of the technical complications is that none of the bundles we define
below are allowed to use the hermitian metric on $V$, orthogonal
complements, orthonormal bases etc., because we wish to remain in the
$\C$-algebraic category. As a general reference for this section and the
next, the reader may consult \cite{Har}.

\bigskip

\begin{definition} For $0\neq v\in V=\eucl{C}{4}$, we will denote the point $[v]\in
\project{C}{3}$ by $v$, whenever no confusion is likely, to simplify 
notation. We have already denoted the vector subspace 
$\spn\{v, Av, A^{\ast}v\}\subset \eucl{C}{4}$ as $W(v)$. Further define 
$W_{3}(v):=W(v)+AW(v)$, and 
$\widetilde{W}_{3}(v):=W(v)+A^{\ast}W(v)$. Clearly both $W_{3}(v)$ and 
$\widetilde{W}_{3}(v)$ contain $W(v)$. 

\bigskip

Since $A$ and $A^{\ast}$ have 
no common eigenvectors, we have $\dim\,W(v)\geq 2$ for all 
$v\in \project{C}{3}$, and $\dim\,W(v)=2$ for all $v\in C$, because of the 
defining equation $v\wedge Av\wedge A^{\ast}v=0$ of $C$. Also, since 
$\dim\,W(v)=2=\dim\,AW(v)$ for $v\in C$, and  
since $0\neq Av\in W(v)\cap AW(v)$, we have $\dim\,W_{3}(v)\leq 3$ for
all $v\in C$. Similarly $\dim\widetilde{W}_{3}(v)\leq 3$ for all $v\in
C$. 
\end{definition}

\bigskip

If there exists a  $v\in C$ such that $\dim\,W_{3}(v)=2$, then we are 
done. For, in this case $W_{3}(v)$ must equal $W(v)$ since it contains 
$W(v)$. Then the dimension 
$\dim\,\widetilde{W}_{3}(v)=2$ or $=3$. If it is 2, $W(v)$ will be a 2-dimensional invariant space for both 
$A$ and $A^{\ast}$, and the main theorem will follow by lemma \ref{invlemma}.
If $\dim\,\widetilde{W}_{3}(v)=3$, then the flag:
$$
0=W_{0}\subset W_{1}={\Bbb C}v \subset W_{2}=W(v)\subset W_{3}=
\widetilde{W}_{3}(v)\subset W_{4}=V 
$$
satisfies the requirements of (ii) in the lemma \ref{flaglemma}, and we are
done. Similarly, if there exists a $v\in C$ with 
$\dim\,\widetilde{W}_{3}(v)=2$, we are again done. Hence we may assume
that:
\begin{eqnarray}
\dim\,W_{3}(v)=\dim\,\widetilde{W}_{3}(v)=3\;\;\;\mbox{for
all}\;\;\;v\in C
\label{dimW3}
\end{eqnarray}

\bigskip

In the light of the above:

\bigskip

\begin{remark} We are reduced to the situation where the following 
condition holds:

\bigskip 
\noin
For each $v\in C$, $\dim\,W(v)=2,\;
\dim\,W_{3}(v)=\dim\,\widetilde{W}_{3}(v)=3$.

\label{dimassump}
\end{remark}

\bigskip

Now our main task is to prove that there exists a $v\in C$ such that the
two 3-dimensional subspaces $W_{3}(v)$ and $\widetilde{W}_{3}(v)$ are the
{\em same}. In that event, the flag 
$$
0=W_{0}\subset W_{1}={\Bbb C}v\subset W_{2}=W(v)\subset W_{3}=W(v)+AW(v)=
W(v)+A^{\ast}W(v)\subset W_{4}=V
$$
will meet the requirements of (ii) of the lemma \ref{flaglemma}. The remainder
of this discussion is aimed at proving this. 

\bigskip

\begin{definition} Denote the trivial rank $4$ algebraic bundle on
$\project{C}{3}$ by $\bun{\project{C}{3}}^{4}$, with fibre $V=\eucl{C}{4}$ 
at each point (following standard
algebraic geometry notation). Similarly, $\bun{C}^{4}$ is the 
trivial bundle on $C$. In $\bun{\project{C}{3}}^{4}$, there is the
tautological line-subbundle $\bun{\project{C}{3}}(-1)$, whose fibre at $v$ is
${\Bbb C}v$. Its restriction to the curve $C$ is denoted as 
${\cal W}_{1}:=\bun{C}(-1)$.

\bigskip

There are also the line subbundles
$A\bun{\project{C}{3}}(-1)$ (respectively $A^{\ast}\bun{\project{C}{3}}(-
1)$) of $\bun{\project{C}{3}}^{4}$, whose fibre at $v$ is $Av$
(respectively $A^{\ast}v$). Both are isomorphic to
$\bun{\project{C}{3}}(-1)$ (via the global linear automorphisms $A$ 
(resp. $A^{\ast}$) of $V$). Similarly, their restrictions $A\bun{C}(-
1)$, $A^{\ast}\bun{C}(-1)$, both isomorphic to $\bun{C}(-1)$. {\em Note that
throughout what follows, bundle isomorphism over any variety $X$ will mean 
algebraic isomorphism, i.e. isomorphism of the corresponding sheaves of 
algebraic sections as $\bun{X}$-modules}. 

\bigskip

Denote the rank 2 algebraic bundle with fibre $W(v)\subset V$ at $v\in C$ 
as  ${\cal W}_{2}$. It is an algebraic subbundle of $\bun{C}^{4}$, for its 
sheaf of sections is the restriction of the subsheaf 
$$
\bun{\project{C}{3}}(-1) + A\bun{\project{C}{3}}(-1)+ 
A^{\ast}\bun{\project{C}{3}}(-1)\subset \bun{\project{C}{3}}^{4}
$$
to the curve $C$, which is precisely the subvariety of $\project{C}{3}$
on which the sheaf above is locally free of rank 2 (=rank 2 algebraic
bundle). 

\bigskip

Denote the rank 3 algebraic 
subbundle of $\bun{C}^{4}$ with fibre $W_{3}(v)=W(v)+AW(v)$ (respectively 
$\widetilde{W}_{3}(v)=W(v)+A^{\ast}W(v)$) by ${\cal W}_{3}$ (respectively 
$\widetilde{{\cal W}}_{3}$). Both ${\cal W}_{3}$ and 
$\widetilde{\cal W}_{3}$ are of rank 3 on $C$ because of 
(\ref{dimassump}) above, and both contain ${\cal W}_{2}$ as a subbundle. 
We denote the line bundles $\ext^{2}{\cal
W}_{2}$ by ${\cal L}_{2}$, and $\ext^{3}{\cal W}_{3}$ (resp.
$\ext^{3}\widetilde{W}_{3}$) by ${\cal L}_{3}$ (resp. 
$\widetilde{\cal L}_{3}$). 
Then ${\cal L}_{2}$ is a line subbundle of $\ext^{2}\bun{C}^{4}$, and ${\cal
L}_{3}, \;\widetilde{\cal L}_{3}$ are line subbundles of 
$\ext^{3}\bun{C}^{4}$.

\bigskip
 
Finally, for $X$ any variety, with a bundle ${\cal E}$ on $X$ which is a 
subbundle of a trivial bundle $\bun{X}^{m}$, the {\em annihilator} 
of ${\cal E}$  is defined as: 
$$
\mbox{Ann}{\cal E}=\{\phi\in \hom_{X}(\bun{X}^{m},\bun{X}):
\phi({\cal E})=0\} 
$$
Clearly, by taking $\hom_{X}(-\,,\bun{X})$ of the exact sequence 
$$
0\rarr {\cal E}\rarr \bun{X}^{m}\rarr \bun{X}^{m}/{\cal E}\rarr 0
$$
the bundle 
$$
\mbox{Ann}{\cal E}\simeq \hom_{X}(\bun{X}^{m}/{\cal E},\bun{X})
=(\bun{X}^{m}/{\cal E})^{\ast}
$$
where $\ast$ always denotes the (complex) dual bundle. 
\label{bunCdef}
\end{definition}

\bigskip

\begin{lemma}
 Denote the bundle ${\cal W}_{3}/{\cal W}_{2}$ (resp. 
$\widetilde{{\cal W}}_{3}/{\cal W}_{2}$) by $\Lambda$ 
(resp. $\widetilde{\Lambda}$). Then we have the following identities of 
bundles on $C$:

\begin{description}
\itm{i} 
\vspace{.2in}
\begin{eqnarray*}
&0&\rarr  {\cal W}_{2} \rarr {\cal W}_{3} \rarr \Lambda \rarr
0\\
\vspace{.2in}
&0&\rarr {\cal W}_{2} \rarr \widetilde{{\cal W}}_{3} \rarr \widetilde{\Lambda}
\rarr 0 \\
\vspace{.2in}
&0&\rarr {\cal L}_{3} \stackrel{i}{\rarr} \mbox{Ann}{\cal W}_{2}
\stackrel{\pi}{\rarr} \Lambda^{\ast} \rarr
0\\
\vspace{.2in}
&0&\rarr\widetilde{{\cal L}}_{3} \stackrel{\til{i}}{\rarr}
\mbox{Ann}{\cal W}_{2}\stackrel{\til{\pi}}{\rarr}\widetilde{\Lambda}^{\ast} 
\rarr 0
\end{eqnarray*}

\bigskip

\itm{ii} 
$$
{\cal L}_{3}\simeq {\cal L}_{2}\otimes
\Lambda\;\;\;\mbox{and}\;\;\;
\widetilde{{\cal L}}_{3}\simeq {\cal L}_{2}\otimes \widetilde{\Lambda}
$$

\bigskip

\itm{iii} 
$$
\ext^{2}\mbox{Ann}{\cal W}_{2}\simeq \ext^{2}{\cal W}_{2}
$$

\bigskip
\itm{iv}
$$
\Lambda \simeq \widetilde{\Lambda}
$$

\bigskip 
\itm{v}
$$
{\cal L}_{2}\simeq \Lambda \otimes {\cal O}_{C}(-1)\simeq
\widetilde{\Lambda}\otimes {\cal O}_{C}(-1)
$$

\itm{vi}
$$
\hom_{C}({\cal L}_{3},\widetilde{\Lambda}^{\ast})\simeq 
{\cal L}_{2}^{\ast}\otimes \widetilde{\Lambda}^{\ast\,2}\simeq 
{\cal L}_{2}^{\ast\,3}\otimes {\cal O}_{C}(-2)
$$

\end{description}
\label{bunlemma}
\end{lemma}

\bigskip
\noin
{\bf Proof:}
From the definition of $\Lambda$, we have the exact sequence:
$$
0 \rarr  {\cal W}_{2} \rarr {\cal W}_{3} \rarr \Lambda \rarr 0
$$
from which it follows that:
$$
0\rarr \Lambda \rarr {\bun{C}}^{4}/{\cal W}_{2}\rarr 
{\bun{C}}^{4}/{\cal W}_{3}\rarr 0
$$
is exact. Taking $\hom_{C}(-,\,\bun{C})$ of this exact sequence yields the 
exact sequence:

$$
0\rarr \mbox{Ann}{\cal W}_{3}\rarr \mbox{Ann}{\cal W}_{2}\rarr
\Lambda^{\ast}\rarr 0
$$

Now, via the canonical isomorphism $\ext^{3}V\rarr V^{\ast}$ which
arises from the non-degenerate pairing 

$$
\ext^{3}V\otimes V\rarr\ext^{4}V\simeq {\Bbb C}
$$
it is clear that $\mbox{Ann}{\cal W}_{3} \simeq \ext^{3}{\cal W}_{3}={\cal L}_{3}$. 

\bigskip

Thus the first and third exact sequences of (i) follow. The 
proofs of the second and fourth are similar. From the 
first exact sequence in (i), it follows that $\ext^{3}{\cal W}_{3}\simeq
\ext^{2}{\cal W}_{2}\otimes \Lambda$. This implies the first identity of
(ii). Similarly the second exact sequence of (i) implies the other identity 
of (ii). 

\bigskip 

Since for every line bundle $\gamma$, $\gamma\otimes \gamma^{\ast}$ is trivial, we 
get from the first identity of (ii) that ${\cal L}_{2}\simeq {\cal
L}_{3}\otimes \Lambda^{\ast}$. From
third exact sequence in (i) it follows that 
$\ext^{2}\mbox{Ann}{\cal W}_{2}\simeq {\cal L}_{3}\otimes \Lambda^{\ast}$, 
and this implies (iii). 

\bigskip

To see (iv), note that 
$$
\Lambda\simeq \frac{{\cal W}_{2}+A{\cal W}_{2}}{{\cal W}_{2}}
\simeq\frac{A{\cal W}_{2}}{A{\cal W}_{2}\cap {\cal W}_{2}}
$$
The automorphism $A^{-1}$ of $V$ makes the last bundle on the right isomorphic
to the line bundle ${\cal W}_{2}/({\cal W}_{2}\cap A^{-1}{\cal W}_{2})$ 
(note all these operations are happening inside the rank 4 trivial bundle
$\bun{C}^{4}$). Similarly, $\til{\Lam}$ is isomorphic (via the global
isomorphism $A^{\ast\,-1}$ of $V$) to the line bundle 
${\cal W}_{2}/({\cal W}_{2}\cap A^{\ast\,-1}{\cal W}_{2})$. But for each 
$v\in C$, 
$W(v)\cap A^{-1}W(v)={\Bbb C}v=W(v)\cap A^{\ast\,-1}W(v)$, 
from which it follows that the line subbundles 
${\cal W}_{2}\cap A^{-1}{\cal W}_{2}$ and 
${\cal W}_{2}\cap A^{\ast\,-1}{\cal W}_{2}$ of ${\cal W}_{2}$ are the same
($={\cal W}_{1}\simeq \bun{C}(-1)$). 
Thus $\Lambda\simeq \widetilde{\Lambda}$, proving (iv). 

\bigskip
To see (v), we need another exact sequence. For each $v\in C$, we noted 
in the proof of (iv) above that ${\Bbb C}v=W(v)\cap 
A^{-1}W(v)$. Thus the sequence
of bundles:
$$
0\rarr {\cal O}_{C}(-1)\rarr {\cal W}_{2}\rarr 
\fr{{\cal W}_{2}}{{\cal W}_{2}\cap A^{-1}{\cal W}_{2}}\rarr 0
$$
is exact. But, as we noted in the proof of (iv) above, the bundle on the
right is isomorphic to $\Lambda$, so that 
$$
0\rarr {\cal O}_{C}(-1)\rarr {\cal W}_{2}\rarr 
\Lambda\rarr 0
$$
is exact. Hence ${\cal L}_{2}=\ext^{2}{\cal W}_{2}\simeq \Lambda\otimes 
{\cal O}_{C}(-1)$. The other identity follows from (iv), thus proving
(v). 

\bigskip

To see (vi) note that we have 
by (ii) ${\cal L}_{3}^{\ast}\simeq {\cal
L}_{2}^{\ast}\otimes \Lambda^{\ast}$. Thus 
$$
\hom_{C}({\cal L}_{3},\widetilde{\Lambda}^{\ast})\simeq 
{\cal L}_{3}^{\ast}\otimes 
\widetilde{\Lambda}^{\ast} \simeq {\cal
L}_{2}^{\ast}\otimes \Lambda^{\ast}\otimes \widetilde{\Lambda}^{\ast}
$$
However, since by (iv), ${\Lambda}\simeq \widetilde{{\Lambda}}$, we have 
$\hom_{C}({\cal L}_{3},\widetilde{\Lambda}^{\ast})\simeq 
{\cal L}_{2}^{\ast}\otimes 
{\Lambda}^{\ast\,2}$. Now, substituting $\Lambda^{\ast}=
{\cal L}_{2}^{\ast}\otimes {\cal O}_{C}(-1)$ from (v), we have the rest
of (vi). Hence the lemma. \hfill$\Box$

\bigskip

We need one more bundle identity: 

\begin{lemma}
There is a bundle isomorphism:
$$
{\cal L}_{2}\simeq g^{\ast}{\cal O}_{D}(1) \otimes {\cal O}_{C}(-2)
$$
\label{bundleiden}
\end{lemma}

\bigskip
\noin
{\bf Proof:} 
When $[t]=[t_{0}:t_{1}:t_{2}] = g([v])$, we saw in (\ref{imageBa}) 
that the linear map 
$B(v,-):\eucl{C}{3}\rarr \eucl{C}{4}$ acquires a 1-dimensional kernel,
which is precisely the line ${\Bbb C}t$, which is the fibre of ${\cal
O}_{D}(-1)$ at $[t]$. The image of $B(v,-)$ was the 2-dimensional 
span $W(v)$ of $v, Av, A^{\ast}v$, as noted there. Thus for $v\in C$, $B(-,-)$ induces a
canonical isomorphism of vector spaces:
$$
{\cal O}_{C}(-1)_{v}\otimes\left(\eucl{C}{3}/{\cal O}_{D}(-1)\right)_{g(v)} \rarr
{\cal W}_{2,v}
$$
which, being defined by the global map $B(-,-)$, gives an isomorphism of
bundles:
$$
{\cal O}_{C}(-1)\otimes g^{\ast}\left({\cal O}_{D}^{3}/{\cal O}_{D}(-1)\right)
\simeq {\cal W}_{2}
$$
From the short exact sequence:
$$
0\rarr {\cal O}_{D}(-1)\rarr {\cal O}_{D}^{3}\rarr 
{\cal O}_{D}^{3}/{\cal O}_{D}(-1)\rarr 0
$$
it follows that $\ext^{2}({\cal O}_{D}^{3}/{\cal O}_{D}(-1))\simeq {\cal
O}_{D}(1)$. Thus:
\begin{eqnarray*}
{\cal L}_{2} &=& \ext^{2}{\cal W}_{2}
\simeq g^{\ast}\left(\ext^{2}({\cal O}_{D}^{3}/{\cal O}_{D}(-1))\right)\otimes {\cal
O}_{C}(-2)\\
&\simeq& g^{\ast}{\cal O}_{D}(1)\otimes {\cal O}_{C}(-2)
\end{eqnarray*}

\bigskip

This proves the lemma. \hfill $\Box$

\bigskip
\section{Degree computations}

In this section, we compute the degrees of the various line bundles
introduced in the previous section. 

\bigskip

\begin{definition}
Note that an {\em irreducible} complex projective curve $C$, as a topological
space, is a canonically oriented pseudomanifold of real dimension 2, and
has a canonical generator $\mu_{C}\in H_{2}(C,{\Bbb Z})={\Bbb Z}$.
Indeed, it is the image  
$\pi_{\ast}\mu_{\til{C}}$, where $\pi:\widetilde{C}\rarr C$ is the
normalisation map, and $\mu_{\widetilde{C}}\in H_{2}(\widetilde{C}, {\Bbb Z})={\Bbb Z}$ 
is the canonical orientation class for the smooth connected 
compact complex manifold $\widetilde{C}$, where 
$\pi_{\ast}:H_{2}(\widetilde{C},{\Bbb Z})\rarr H_{2}(C,{\Bbb Z})$ is an 
isomorphism for elementary topological reasons. 

\bigskip

If $C=\cup_{\alp=1}^{r}C_{\alp}$ is a projective curve of pure dimension 1, 
with the curves $C_{\alp}$ as irreducible components, then since the
intersections $C_{\alp}\cap C_{\beta}$ are finite sets of points (or empty), 
$H_{2}(C,{\Bbb Z})=\oplus_{\alp}H_{2}(C_{\alp}, {\Bbb Z})$. Letting 
$\mu_{\alp}$ denote the canonical orientation classes of $C_{\alp}$ as above, 
there is a {\em unique class} 
$\mu_{C}=\sum_{\alp}\mu_{\alp} \in H_{2}(C,{\Bbb Z})$. Thinking of $C$
as an oriented 2-pseudomanifold, $\mu_{C}$ is just the sum of all the
oriented 2-simplices of $C$.

\bigskip

If ${\cal F}$ is
a complex line bundle on $C$, it has a first Chern class 
$c_{1}({\cal F})\in H^{2}(X, {\Bbb Z})$, and the 
{\em degree} of ${\cal F}$ is defined by:
$$
\mbox{deg}\,{\cal F}=\innp{c_{1}({\cal F})}{\mu_{C}}\in {\Bbb Z}
$$
It is known that a complex line bundle on a pseudomanifold is 
topologically trivial iff its first Chern class is zero. In particular,
if an algebraic line bundle on a projective variety has non-zero
degree, then it is topologically (and hence algebraically) non-trivial. 

\bigskip

Finally, if $i:C\hookrightarrow \project{C}{n}$ is an (algebraic) embedding 
of a curve in some projective space, we define the degree of the bundle ${\cal
O}_{C}(1)=i^{\ast}{\cal O}_{\project{C}{n}}(1)$ as the {\em degree of the 
curve} $C$ (in $\project{C}{n}$). We note that $[C]:=i_{\ast}(\mu_{C})\in H_{2}(\project{C}{n},{\Bbb Z})$
is called the {\em fundamental class} of $C$ in $\project{C}{n}$, and by
definition $\mbox{deg}\,C=\innp{c_{1}({\cal O}_{C}(1))}{\mu_{C}}=
\innp{c_{1}({\cal O}_{\project{C}{n}}(1))}{[C]}$. Geometrically, one  
intersects $C$ with a generic hyperplane, which intersects $C$ away from 
its singular locus in a finite set of points, and then counts these points 
of intersection with their multiplicity. 

\bigskip

More generally, a complex projective variety $X\subset \project{C}{n}$ 
of complex dimension $m$ has a unique orientation class $\mu_{X}\in 
H_{2m}(X,{\Bbb Z})$. Its image in
$H_{2m}(\project{C}{n},{\Bbb Z})$ is denoted $[X]$, and the degree $\mbox{deg}\,X$ 
of $X$ is defined as $\innp{(c_{1}({\cal O}_{\project{C}{n}}(1)))^{m}}{[X]}$. It
is known that if $X=V(F)$ for a homogeneous polynomial $F$ of degree $d$,
then $\mbox{deg}\,X=d$.

\label{degreedef}
\end{definition}

\bigskip

We need the following remark later on:

\begin{remark} If $f:C\rarr D$ is a regular isomorphism of complex 
projective curves $C$ and $D$, both of pure dimension 1, and 
if ${\cal F}$ is a complex line
bundle on $D$, then $\mbox{deg}\,f^{\ast}{\cal F}=\mbox{deg}\,{\cal F}$. 
This is because $f_{\ast}(\mu_{C})=\mu_{D}$, so that 
$$
\mbox{deg}\,{\cal F}=\innp{c_{1}({\cal F})}{\mu_{D}}=
\innp{c_{1}({\cal F})}{f_{\ast}\mu_{C}}=
\innp{f^{\ast}c_{1}({\cal F})}{\mu_{C}}=
\innp{c_{1}(f^{\ast}{\cal F})}{\mu_{C}}=\mbox{deg}\,f^{\ast}{\cal F}
$$
\label{degremark}
\end{remark}

Now we can compute the degrees of all the line bundles introduced.

\newpage
\begin{lemma}The degrees of the various line bundles above are as
follows:

\begin{description}
\itm{i} 
$\mbox{deg}\,{\cal O}_{C}(1)=\mbox{deg}\,C = 6$

\bigskip

\itm{ii} 
$\mbox{deg}\,{\cal O}_{D}(1)=\mbox{deg}\,D = 4$

\bigskip

\itm{iii} 
$\mbox{deg}\,{\cal L}_{2}^{\ast}= 8$

\bigskip

\itm{iv} $\mbox{deg}\,\hom_{C}({\cal L}_{3},
\widetilde{\Lambda}^{\ast})=\mbox{deg}\left({\cal L}_{2}^{\ast 3}\otimes
{\cal O}_{C}(-2)\right)=12$
\end{description}

\label{degreelemma}
\end{lemma}

\bigskip
\noin 
{\bf Proof:}

We denote the image of orientation class $\mu_{\Gamma}$ of 
the curve $\Gamma$ (see
definition \ref{defgamma} for the definition of $\Gamma$)
in $H_{2}(\project{C}{3}\times \project{C}{2}, {\Bbb Z})$
by $[\Gamma]$. By the part (iv) of the Lemma \ref{poslemma2}, we have
that the homology class $[\Gamma]$ is the same as the homology class of
the intersection cycle defined by the four divisors $D_{i}:=(B_{i}(v,t)=0)$
inside $H_{2}(\project{C}{3}\times \project{C}{2}, {\Bbb Z})$. By the
generalised Bezout theorem in $\project{C}{3}\times\project{C}{2}$, the
homology class of the last-mentioned intersection cycle is the homology
class Poincare-dual to the cup product 
$$
d:=d_{1}\cup d_{2}\cup d_{3}\cup
d_{4}
$$
where $d_{i}$ is the first Chern class of the the line bundle
$L_{i}$ corresponding to $D_{i}$, for $i=1,2,3,4$. (See \cite{Sha}, p.
237, Ex.2). 

\bigskip

Since each $B_{i}(v,t)$ is separately linear in $v$, $t$, the line
bundle defined by the divisor $D_{i}$ 
is the bundle $\pi_{1}^{\ast}{\cal O}_{\project{C}{3}}(1)\otimes 
\pi_{2}^{\ast}{\cal O}_{\project{C}{2}}(1)$,
where $\pi_{1},\;\pi_{2}$ are the projections to $\project{C}{3}$ and 
$\project{C}{2}$ respectively. If we denote the hyperplane classes which 
are the generators of the cohomologies 
$H^{2}(\project{C}{3},{\Bbb
Z})$ and $H^{2}(\project{C}{2},{\Bbb Z})$ by $x$ and $y$ respectively,
we have:
$$
d_{i}=c_{1}(L_{i})=\pi_{1}^{\ast}(x) + \pi_{2}^{\ast}(y)
$$

\bigskip

Then we have, from the cohomology ring structures of $\project{C}{3}$
and $\project{C}{2}$ that $x\cup x\cup x\cup x =y\cup y\cup y=0$. Hence
the cohomology class in $H^{8}(\project{C}{3}\times\project{C}{2},{\Bbb
Z})$ given by the cup-product of $d_{i}$ is: 
$$
d:=d_{1}\cup d_{2}\cup d_{3}\cup d_{4}= (\pi_{1}^{\ast}(x) +
\pi_{2}^{\ast}(y))^{4} = 
4\pi_{1}^{\ast}(x^{3})\pi_{2}^{\ast}(y) +
6\pi_{1}^{\ast}(x^{2})\pi_{2}^{\ast}(y^{2})
$$
where $x^{3}=x\cup x\cup x$ etc. By part (ii) of the lemma \ref{poslemma2},
 the map $\pi_{1}:\Gamma\rarr C$ is an isomorphism, so 
applying the remark \ref{degremark} to it, we have:

\begin{eqnarray}
\mbox{deg}\,{\cal O}_{C}(1) &=& \mbox{deg}\,\pi_{1}^{\ast}{\cal O}_{C}(1)\nonumber\\
&=&\innp{c_{1}(\pi_{1}^{\ast}({\cal O}_{\project{C}{3}}(1))}{[\Gamma]}
\nonumber\\
&=&\innp{c_{1}(\pi_{1}^{\ast}({\cal O}_{\project{C}{3}}(1))\cup d}
{[\project{C}{3}\times\project{C}{2}]}\nonumber\\
&=&\innp{\pi_{1}^{\ast}(x)\cup 
\left(4\pi_{1}^{\ast}(x^{3})\pi_{2}^{\ast}(y) +
6\pi_{1}^{\ast}(x^{2})\pi_{2}^{\ast}(y^{2})\right)}{[\project{C}{3}\times
\project{C}{2}]}\nonumber\\
&=&\innp{6\pi_{1}^{\ast}(x^{3})\cup\pi_{2}^{\ast}(y^{2})}{[\project{C}{3}
\times\project{C}{2}]}\nonumber\\
&=& 6
\label{cupproducts}
\end{eqnarray}
where we have used the Poincare duality cap-product relation 
$[\Gamma]= [\project{C}{3}\times\project{C}{2}]\cap d$ mentioned above,
and that $\pi_{1}^{\ast}(x^{3})\cup\pi_{2}^{\ast}(y^{2})$ is the generator of
$H^{10}(\project{C}{3}\times\project{C}{2}, {\Bbb Z})$, so evaluates to 
1 on the orientation class $[\project{C}{3}\times\project{C}{2}]$, and 
$x^{4}=0$.  This proves (i).

\bigskip

The proof of (ii) is similar, we just replace $C$ by $D$, and $\pi_{1}$
by $\pi_{2}$, and $\pi_{1}^{\ast}(x)$ by $\pi_{2}^{\ast}(y)$ in the 
equalities of (\ref{cupproducts}) above, and get 4 
(as one should expect, since $D$ is defined by a degree 4 homogeneous
polynomial in $\project{C}{2}$). This proves (ii).

\bigskip

For (iii), we use the identity of lemma \ref{bundleiden} that ${\cal
L}_{2}=g^{\ast}{\cal O}_{D}(1)\otimes {\cal O}_{C}(-2)$, and the remark 
\ref{degremark} applied to the isomorphism of curves $g:C\rarr D$ (part
(iii) of the lemma \ref{poslemma2}) to conclude that $\mbox{deg}\,{\cal
L}_{2}=\mbox{deg}\,D-2\mbox{deg}\,C= 4-12=-8$, by (i) and (ii) above, 
so that $\mbox{deg}\,{\cal L}_{2}^{\ast}=8$. 

\bigskip 

For (iv), we have by (vi) of the lemma \ref{bunlemma} that 
$\hom_{C}({\cal L}_{3},\widetilde{\Lambda}^{\ast})\simeq {\cal
L}_{2}^{\ast 3}\otimes {\cal O}_{C}(-2)$, so that its degree is 
$3\mbox{deg}\,{\cal L}_{2}^{\ast}-2\mbox{deg}\,C=24-12 =12$ by 
(i) and (iii) above. 

This proves the lemma. \hfill $\Box$

\bigskip

From (iv) of the lemma above, we have the:

\bigskip

\begin{corollary}
The line bundle $\hom_{C}({\cal L}_{3}, \widetilde{\Lambda}^{\ast})$ is 
a non-trivial line bundle.
\label{poslemma3}
\end{corollary}

\bigskip

\section{Proof of the Main Theorem}

\bigskip
\noin
{\bf Proof of Theorem \ref{mainthm}:}  By the third exact sequence in (i) of the lemma
\ref{bunlemma}, we have a bundle morphism $s$ of line bundles on $C$ defined as 
the composite:
$$
\mbox{Ann}{\cal W}_{3}={\cal L}_{3} \stackrel{i}{\rarr} \mbox{Ann}{\cal W}_{2}
\stackrel{\widetilde{\pi}}{\rarr}\widetilde{\Lambda}^{\ast}=\mbox{Ann}{\cal
W}_{2}/\mbox{Ann}\widetilde{{\cal W}}_{3}
$$
which vanishes at $v\in C$ if and only if the fibre $\mbox{Ann}{\cal W}_{3\,,v}$ is
equal to the fibre $\mbox{Ann}\widetilde{\cal W}_{3\,,v}$ inside
$\mbox{Ann}{\cal W}_{2\,,v}$. At such a point $v\in C$, we have 
$\mbox{Ann}{\cal W}_{3,v}=\mbox{Ann}\widetilde{\cal W}_{3,v}$, so that 
${\cal W}_{3,v}=W(v)+AW(v)=\widetilde{\cal W}_{3,v}=W(v)+A^{\ast}W(v)$. But
this morphism $s$ is a global section of the bundle 
$\hom_{C}({\cal L}_{3},\widetilde{\Lambda}^{\ast})$, which is not a
trivial bundle by the corollary \ref{poslemma3} of the last section. 

\bigskip

Thus there exists a $v\in C$, satisfying $s(v)=0$, and consequently the flag
\begin{eqnarray*}
0\subset W_{1}={\cal W}_{1\,,v}&=&{\Bbb C}v\subset  W_{2}={\cal
W}_{2\,,v}=W(v)=\spn\{v,Av,A^{\ast}v\}\\
\subset W_{3}&=&{\cal W}_{3\,,v}=W(v)+AW(v)=W(v)+A^{\ast}W(v)=
\widetilde{\cal W}_{3,v}\subset
W_{4}=V=\eucl{C}{4}
\end{eqnarray*}
satisfies the requirements of (ii) of lemma \ref{flaglemma}, and the main theorem
\ref{mainthm} follows. \hfill $\Box$

\bigskip
\begin{remark} Note that since $\dim\,C=1$, then the set of points $v\in C$ 
such $s(v)=0$, where $s$ is the section above, will be a finite set. Then  
the set of flags that satisfy (ii) of lemma \ref{flaglemma} which 
tridiagonalise $A$ of the kind considered above (viz. $A$ satisfying the
assumptions of \ref{assump}), will only be finitely many (at most 12 in number!). 
\end{remark}

\bigskip
{\bf Acknowledgments} I am grateful to Bhaskar Bagchi for
posing the problem to me, and to B.V. Rajarama Bhat and J. Holbrook 
for pointing me to the relevant literature. I am deeply grateful to 
the referee, whose valuable comments have led to the elimination of
grave errors, and a substantial streamlining of this paper.

\bigskip
\noin

\end{document}